\documentclass[dvips,ejs]{imsart}

\usepackage{amsthm,amsmath,amsfonts,amssymb}
\usepackage[french]{babel}
\usepackage[latin1]{inputenc}
\RequirePackage[colorlinks,citecolor=blue,urlcolor=blue]{hyperref}

\RequirePackage{hypernat}

\doi{10.1214/154957804100000000}
\pubyear{0000}
\volume{0}
\firstpage{0}
\lastpage{0}

\startlocaldefs
\endlocaldefs

\begin{document}
\newcommand{\dd}{\delta}
\newcommand{\e}{\epsilon}
\newcommand{\alp}{\alpha}
\newcommand{\lab}{\lambda}
\newcommand{\gam}{\gamma}
\newcommand{\sig}{\sigma}
\newcommand{\tht}{\theta}
\newcommand{\Lab}{\Lambda}
\newcommand{\Gam}{\Gamma}
\newcommand{\Sig}{\Sigma}
\newcommand{\T}{\Theta}


\newcommand{\sli}{\sum\limits}
\newcommand{\sliin}{\sum\limits_{i=1}^n}
\newcommand{\sliid}{\sum\limits_{i=1}^d}
\newcommand{\sliik}{\sum\limits_{i=1}^k}
\newcommand{\sliiN}{\sum\limits_{i=1}^N}
\newcommand{\slijk}{\sum\limits_{j=1}^k}
\newcommand{\proliin}{\prod\limits_{i=1}^n}
\newcommand{\proliid}{\prod\limits_{i=1}^d}
\newcommand{\proliik}{\prod\limits_{i=1}^k}
\newcommand{\proliiN}{\prod\limits_{i=1}^N}
\newcommand{\prolijk}{\prod\limits_{j=1}^k}
\newcommand{\ili}{\int}
\newcommand{\proli}{\prod\limits}
\newcommand{\bculi}{\bigcup\limits}
\newcommand{\bcali}{\bigcap\limits}


\newcommand{\ls}{\limsup}
\newcommand{\li}{\liminf}
\newcommand{\limn}{\lim_{n\rightarrow\infty}\;}
\newcommand{\limk}{\lim_{k\rightarrow\infty}\;}
\newcommand{\lsn}{\limsup_{n\rightarrow\infty}}
\newcommand{\lin}{\liminf_{n\rightarrow\infty}}
\newcommand{\lsk}{\limsup_{k\rightarrow\infty}}
\newcommand{\lik}{\liminf_{k\rightarrow\infty}}
\newcommand{\rar}{\rightarrow}
\newcommand{\lar}{\leftarrow}
\newcommand{\cvps}{\rightarrow_{p.s.}\;}
\newcommand{\cvpr}{\rightarrow_{P}\;}
\newcommand{\cvloi}{\rightarrow_{\mcal{L}}\;}
\newcommand{\kif}{k\rightarrow\infty}
\newcommand{\nif}{n\rightarrow\infty}

\newcommand{\coe}{[}
\newcommand{\cfe}{]}
\newcommand{\aoo}{\Big\{}
\newcommand{\aff}{\Big\}}
\newcommand{\coo}{\Big [}
\newcommand{\cff}{\Big]}
\newcommand{\poo}{\Big (}
\newcommand{\pff}{\Big)}
\newcommand{\po}{\big (}
\newcommand{\pf}{\big)}
\newcommand{\ao}{\big \{}
\newcommand{\af}{\big \}}
\newcommand{\co}{\big [}
\newcommand{\cf}{\big ]}
\newcommand{\pooo}{\bigg (}
\newcommand{\pfff}{\bigg)}
\newcommand{\aooo}{\bigg \{}
\newcommand{\afff}{\bigg \}}
\newcommand{\cooo}{\bigg [}
\newcommand{\cfff}{\bigg ]}
\newcommand{\poooo}{\Bigg (}
\newcommand{\pffff}{\Bigg)}
\newcommand{\aoooo}{\Bigg \{}
\newcommand{\affff}{\Bigg \}}
\newcommand{\coooo}{\Bigg [}
\newcommand{\cffff}{\Bigg ]}
\newcommand{\mmi}{\mid\mid}
\newcommand{\mmmi}{\mid\mid\mid}
\newcommand{\mmii}{\mid\mid_\infty}
\newcommand{\gMid}{\Bigg |}
\newcommand{\Mid}{\Big |}
\newcommand{\Mmi}{\Big|\Big|}
\newcommand{\Mmii}{\Big|\Big|_\infty}


\newcommand{\nk}{n_k}
\newcommand{\nkm}{n_{k-1}}
\newcommand{\ovnk}{\overline{n}_k}
\newcommand{\ovnkm}{\overline{n}_{k-1}}
\newcommand{\nkk}{n_{k+1}}
\newcommand{\hnk}{h_{n_k}}
\newcommand{\bnk}{b_{n_k}}
\newcommand{\bnkk}{b_{n_{k+1}}}


\newcommand{\EEE}{\mathbb{E}}
\newcommand{\NNN}{\mathbb{N}}
\newcommand{\PPP}{ \mathbb{P}}
\newcommand{\CCC}{\mathbb{C}}
\newcommand{\KKK}{\mathbb{K}}
\newcommand{\RRR}{\mathbb{R}}
\newcommand{\UUU}{\mathbb{U}}
\newcommand{\QQQ}{\mathbb{Q}}
\newcommand{\cc}{\mathfrak{c}}
\newcommand{\AAA}{\mathcal{A}}
\newcommand{\BB}{\mathcal{B}}
\newcommand{\FF}{\mathcal{F}}
\newcommand{\TT}{\mathcal{T}}
\newcommand{\GG}{\mathcal{G}}
\newcommand{\BBGG}{\mathcal{B}(\mathcal{G})}
\newcommand{\CC}{\mathcal{C}}
\newcommand{\KK}{\mathcal{K}}
\newcommand{\SSS}{\mathcal{S}}
\newcommand{\PP}{\mathcal{P}}
\newcommand{\HH}{\mathcal{H}}
\newcommand{\NN}{\mathcal{N}}
\newcommand{\MM}{\mathcal{M}}
\newcommand{\DD}{\mathcal{D}}
\newcommand{\LL}{\mathcal{L}}
\newcommand{\VV}{\mathcal{V}}

\newcommand{\ovg}{\overline{g}}
\newcommand{\ovh}{\overline{h}}
\newcommand{\ovE}{\overline{E}}
\newcommand{\ovH}{\overline{H}}
\newcommand{\ovI}{\overline{I}}
\newcommand{\ovJ}{\overline{J}}
\newcommand{\ovK}{\overline{K}}
\newcommand{\ovR}{\overline{R}}
\newcommand{\ovFF}{\overline{\mathcal{F}}}
\newcommand{\ovPPP}{\overline{\mathbb{P}}}

\newcommand{\wt}{\widetilde}
\newcommand{\wtG}{\widetilde{G}}
\newcommand{\wtf}{\widetilde{f}}
\newcommand{\wth}{\mathfrak{h}}
\newcommand{\wtn}{\widetilde{n}}
\newcommand{\wtv}{\widetilde{v}}
\newcommand{\wtA}{\widetilde{A}}
\newcommand{\wtC}{\widetilde{C}}
\newcommand{\wtE}{\widetilde{E}}
\newcommand{\wtF}{\widetilde{F}}
\newcommand{\wtI}{\widetilde{I}}
\newcommand{\wtK}{\widetilde{K}}
\newcommand{\wtJ}{\widetilde{J}}
\newcommand{\wtN}{\widetilde{N}}
\newcommand{\wtP}{\widetilde{P}}
\newcommand{\wtU}{\widetilde{U}}
\newcommand{\wtY}{\widetilde{Y}}
\newcommand{\wtFF}{\widetilde{\mathcal{F}}}
\newcommand{\wtPPP}{\widetilde{\mathbb{P}}}

\newcommand{\ov}{\overline}
\newcommand{\hk}{h_{n_k}}
\newcommand{\lb}{\newline}
\newcommand{\indep}{{\bot}\kern-0.9em{\bot}} 
\newcommand{\beq}{\begin{equation} }
\newcommand{\eeq}{\end{equation} }
\newcommand{\mcal}{\mathcal}
\newcommand{\sq}{\sqrt}
\newcommand{\srl}{\stackrel}
\newcommand{\wap}{(\Omega,\mathcal{A},\rm I\kern-2pt P)} 
\newcommand{\vk}{\vskip10pt}
\newcommand{\norm}{\mid\mid \cdot \mid\mid}
\newcommand{\nono}{\nonumber}
\newcommand{\Cov}{\mathrm{Cov}}
\newcommand{\Var}{\mathrm{Var}}
\newcommand{\Idd}{[0,1)^d}
\newcommand{\suite}{_{n\geq 1}}
\newcommand{\Dnc}{\Delta_{n,\mathfrak{c}}}
\newcommand{\DPnc}{\Delta\Pi_{n,\mathfrak{c}}}
\newcommand{\DPnkc}{\Delta\Pi_{n_k,\mathfrak{c}}}
\newcommand{\mh}{{\mathfrak{h}}}
\newcommand{\mhY}{{\mathfrak{h}_Y}}
\newcommand{\mhYY}{{\mathfrak{h}_{\mid Y \mid_k}}}
\newcommand{\JhY}{J_{\mathfrak{h}_{Y}}}
\newcommand{\GhY}{\Gamma_{J_{\mathfrak{h}_Y}}}
\newcommand{\ii}{\mathbf{i}}
\newcommand{\jj}{\mathbf{j}}
\newcommand{\ppd}{\{1,\ldots,2^p\}^d}
\newtheorem{theo}{Theorem}
\newtheorem{ptheo}{Preuve du théorème}[section]
\newtheorem{lem}{Lemma}[section]
\newtheorem{plem}{Preuve du lemme}[section]
\newtheorem{prop}{Propriété}[section]
\newtheorem{preuveprop}{Preuve de la propriété}[section]
\newtheorem{defi}{Definition}[section]
\newtheorem{propo}{Proposition}[section]
\newtheorem{popo}{Preuve de la proposition}[section]
\newtheorem{coro}{Corollary}[theo]
\newtheorem{pcoro}{Preuve du corollaire}[theo]
\newtheorem{ineg}{Inégalité}[section]
\newtheorem{pineg}{Preuve de l'inégalité}[section]
\newtheorem{rem}{Remark}[section]
\newtheorem{fact}{Fact}[section]
\numberwithin{equation}{section}
\begin{frontmatter}

\title{Non standard functional limit laws for the increments of the compound empirical distribution function}
\runtitle{Non standard behaviour of the compound empirical increments}


\author{\fnms{Myriam} \snm{Maumy}\ead[label=e1]{mmaumy@math.u-strasbg.fr}}
\address{\printead{e1}}
\and
\author{\fnms{Davit} \snm{Varron}\ead[label=e2]{dvarron@univ-fcomte.fr}}
\address{\printead{e2}}

\runauthor{M. Maumy and D. Varron}

\begin{abstract}
Let $(Y_i,Z_i)_{i\geq 1}$ be a sequence of independent,
identically distributed (i.i.d.) random vectors taking values in
$\RRR^k\times\RRR^d$, for some integers $k$ and $d$. Given $z\in
\RRR^d$, we provide a nonstandard functional limit law for the
sequence of functional increments of the compound empirical
process, namely
$$\mathbf{\Delta}_{n,\cc}(h_n,z,\cdot):= \frac{1}{nh_n}\sliin
1_{[0,\cdot)}\poo \frac{Z_i-z}{{h_n}^{1/d}}\pff Y_i.$$ Provided
that $nh_n\sim c\log n $ as $\nif$, we obtain, under some natural
conditions on the conditional exponential moments of $Y\mid Z=z$,
that
$$\mathbf{\Delta}_{n,\cc}(h_n,z,\cdot)\leadsto \Gam\text{  almost surely},$$
where $\leadsto$ denotes the clustering process under the sup norm
on $\Idd$. Here, $\Gam$ is a compact set that is related to the
large deviations of certain compound Poisson processes. 
\end{abstract}



\received{\smonth{1} \syear{0000}}

\tableofcontents

\end{frontmatter}
\section{Introduction and statement of the results}
Let $(Y_i,Z_i)_{i\geq 1}$ be a sequence of independent,
identically distributed (i.i.d.) random vectors taking values in
$\RRR^k\times\RRR^d$, for some integers $k$ and $d$. Given $s,t\in
\RRR^d$ with respective coordinates $s_1,\ldots,s_d$ and
$t_1,\ldots,t_d$, we shall write
$[s,t]:=[s_1,t_1]\times\ldots\times[s_d,t_d]$,
$[s,t):=[s_1,t_1)\times \ldots\times[s_d,t_d)$ and given $a\in
\overline{\RRR}$ we set $[a,t]:=
[a,t_1]\times\ldots\times[a,t_d]$. For each integer $n\geq 1$,
define the compound empirical distribution function as: \beq
\UUU_{n,\cc}(s):=\frac{1}{n}\sliin
 1_{(-\infty,s]}(Z_i)Y_i,\;s\in \RRR^d.\label{UUUnc}\eeq
Here the letter $\cc$ stands for "compound". In this paper, we are
concerned with the asymptotic behaviour of the functional
increments of $\UUU_{n,\cc}$, namely, for fixed $h>0$ and $z\in
\RRR^d,$ \beq \mathbf{\Delta}_{n,\cc}(h,z,s):= \frac{1}{nh}\sliin
1_{[0,s)}\poo \frac{Z_i-z}{h^{1/d}}\pff Y_i,\;s\in \Idd.\eeq Note
that, in the particular case where $k=1$ and $Y_1\equiv 1$, the
$\mathbf{\Delta}_{n,\cc}(h,z,\cdot)$ are no more than the
functional increments of the empirical distribution function,
which have been intensively investigated in the literature (see,
e.g., \cite{MSW,DeheuvelsM2,DeheuvelsM3,Mason1}). Among these investigations, Deheuvels and Mason
(\cite{DeheuvelsM2,DeheuvelsM3}) have established
\textit{nonstandard} functional limit laws for the
$\mathbf{\Delta}_{n,\cc}(h,z,\cdot)$ when $k=1$, $d=1$, $Y_1\equiv
1$ and $Z_1$ is uniformly distributed on $[0,1)$. To cite their
results, we need to introduce some further notations. We shall
write $$\mmi s\mmi_k:= \max\{\mid s_1 \mid,\ldots,\mid s_k \mid
\}$$ for $s\in\RRR^k$, and we define $B_k(\Idd)$ as the space of
all mappings from $\Idd$ to $\RRR^k$ that are bounded. We shall
endow $B_k(\Idd)$ with the usual sup-norm, namely $\mid g\mid_k
:=\sup_{s\in \Idd}\mid g(s)\mid_k.$ Given a convex real function
$\mathfrak{h}$ on $\RRR^k$, we define the following functional on
$B_k(\Idd)$: whenever a function $g$ satisfies $g(0)=0$ and admits
a derivative $g'$ with respect to the Lebesgue measure, set \beq
J_{\mathfrak{h}}(g):=\ili_{\Idd}
\mathfrak{h}(g'(s))ds\label{JZ},\eeq and set
$J_{\mathfrak{h}}(g):=\infty$ if it is not the case. We also
write, for any $c>0$, \beq \Gam_{\mathfrak{h}}(c):=\ao g\in
B_k(\Idd),\;J_{\mathfrak{h}}(g)\le c\label{GamZc}\af .\eeq Now
define the following (Chernoff) function on $[0,\infty)$ :
\beq \mathfrak{h}_1(x):=\left\{%
\begin{array}{ll}
    x\log x -x+1, & \hbox{for $x>0$;} \\
    1, & \hbox{for $x=0$;} \\
    \infty, & \hbox{for $x<0$.} \\
\end{array}%
\right.\label{h}\eeq A sequence $(f_n)$ in a metric space
$(E,\rho)$, is said to be \textit{relatively compact} with limit
set equal to $K$ when $K$ is (non void) compact and the following
assertions are true
\begin{align}
  &\; \limn \inf_{f\in K}d(f_n,f)=0,\label{i}\\
   &\;\forall\;f\in K,\; \lin d(f_n,f)=0\label{ii}.
\end{align}
We shall write this property $x_n\leadsto K$.\lb Throughout this
article, we shall consider a sequence of constants $(h_n)\suite$
satisfying the so called local nonstandard conditions, namely, as
$\nif$,
\begin{tabbing} $\hskip5pt$ \= (HV) $\hskip5pt$ \=
$0<h_n<1,\;h_n\downarrow 0,\;nh_n\uparrow \infty,\; nh_n/\log_2
n\rar c \in (0,\infty).$
\end{tabbing}
Here we have set $\log_2 n:=\log(\log(n\vee 3))$, with the
notation $a\vee b:=\max\{a,b\}$. In a pioneering work, Deheuvels
and Mason \cite{DeheuvelsM3} established a \textit{nonstandard}
functional law of the iterated logarithm for a single functional
increment of the empirical distribution function. With the
notation of the present paper, their theorem can be stated as
follows.\lb
 \begin{fact} [Deheuvels and Mason, 1990]
Let $(h_n)\suite$ be a sequence satisfying (HV) for some constant
$c>0$. Assume that $k=1$, $d=1$, $Y_1\equiv 1$, and that $Z_1$ is
uniformly distributed on $[0,1)$. Then, given $z\in [0,1)$, we
have almost surely
\begin{align}\nono
\mathbf{\Delta}_{n,\cc}(h_n,z,\cdot)\;\leadsto
\Gam_{\mathfrak{h}_1}(1/c).
\end{align}
\end{fact}
Later, Deheuvels and Mason \cite{DeheuvelsM6} extended the just
mentioned result to a more general setting, where $d>1$ and with
fewer assumptions on the law of the $Z_i$, considering the
$\mathbf{\Delta}_{n,\cc}(h_n,z,\cdot)$ as random measures indexed
by a class of sets. The aim of the present paper is to extend the
above mentioned results to the case where the random vectors $Y_i$
are not constant, but do satisfy some assumptions on their
conditional exponential moments given $Z=z$. From now on
$<\cdot,\cdot>$ will always denote the Euclidian scalar product on
$\RRR^k$ and $\lab$ stands for the Lebesgue measure. Define $\CC$
as the class of each $C\subset \RRR^d$ which is the union of $d$
hypercubes of $\RRR^d$, and with $\lab(C)>0$. The two key
assumptions that we shall make upon the law of $(Y_1,Z_1)$ are
stated as follows.
\begin{tabbing}
 $ \hskip 5pt $ \=  (HL1)  $ \hskip 5pt $ \= There exists a constant $f(z)>0$ satisfying, for
 each $C\in \CC$ \\
\>\>  $\displaystyle{\lim_{h\rar 0}  h^{-1}\PPP\poo Z_1\in
z+h^{1/d}C\pff=\lab(C)f(z)}$.\\
 \>  (HL2)  \> There exist two mappings $\LL_Y: \RRR^k\mapsto [0,\infty)$ and\\
 \>\> $\LL_{\mid Y \mid_k}: \RRR\mapsto [0,\infty)$ such that, for each $t\in \RRR^k$  and $t'\in \RRR$ and $C\in \CC$, we have\\
\>\>$\displaystyle{\lim_{h\rar 0}}\EEE\poo \exp\po<t,Y_1>\pf\Mid
Z_1\in z+h^{1/d} C
\pff=\LL_Y(t)$,\\
\>\>$\displaystyle{\mathop{\lim}_{h\rar 0} \EEE\poo \exp\po t'\mid
Y_1\mid_k \pf\Mid Z_1\in z+h^{1/d} C \pff= \LL_{\mid
Y\mid_k}(t')}$.
\end{tabbing}
\begin{rem}
Assumptions (HL1) and (HL2) seem to be the weakest that we can
afford in this context, in regard to the methods we make use of in
this paper. Note that $(HL2)$ implies that $\LL_Y$ is infinitely
differentiable on $\RRR^k$. Some straightforward analysis shows
that these assumptions are fulfilled when the $Y_i$ are bounded by
a constant and when $Z_i$ admit a (version of) density $f$ which
is continuous at $z$. Another interesting case where (HL1) and
(HL2) are fulfilled is a general semi parametric setting which
appears in the following proposition.
\end{rem}
\begin{propo}
Assume that there exist a $\sig$-finite measure $v$ on $\RRR^k$
and an application $f_Y$ from $\RRR^k\times \RRR^d$ to
$[0,\infty)$, such that
\begin{enumerate}
  \item There exists a neighborhood $\VV$ of $z$ such that, for each $z'\in \VV$, the law of $Y\mid Z'=z$ is dominated
  by $v$, with density $f_Y(\cdot,z')$.
  \item For $v$-almost all $y$, the function $z\rar f_Y(y,z)$ is
  continuous on $\VV$.
  \item For each $z'\in \VV$, and each $t\in \RRR^k$ we have
  $$\ili_{\RRR^k} \exp\po <t,u>\pf f_Y(u,z') dv(u)<\infty.$$
  \item $Z$ has a version of density (with respect to the Lebesgue
  measure $\lab$) which is continuous on $\VV$.
\end{enumerate}
Then the random vector $(Y,Z)$ fulfills $(HL1)$ and $(HL2)$.
\end{propo}
\textbf{Proof}: The proof is a straightforward application of
Sch\'eff\'e's lemma. $\Box$.
\begin{rem}
Roughly speaking, assumption (HL2) imposes that the Laplace
transform of the law $Y\mid Z=z$ is finite on $\RRR^k$. One could
argue that this assumption could be weakened. However, it seems
that, when this assumption is dropped, Theorem \ref{T1} (see
below) does not hold anymore under the strong norm $\norm_k$. A
close look at the works of Deheuvels \cite{Deheuvels4} and
Borovkov \cite{Borovkov1} on the functional increments of random
walks leads to the conjecture that the appropriate topology when
$\LL$ is finite only on a neighborhood of $0$  seems to be the
usually called \textit{weak star} topology (see, e.g.
\cite{Deheuvels4}). This topic is however beyond the scope of this
article, and shall be investigated in future works.
\end{rem}
Notice that $\LL_Y$ and $\LL_{\mid Y\mid_k}$ are positive convex
functions when they exist. We now introduce $\mathfrak{h}_Y$
(resp. $\mathfrak{h}_{\mid Y\mid_k}$ ), which is defined as the
Legendre transform of $\LL_Y-1$ (resp. $\LL_{\mid Y\mid_k }-1)$,
namely :
\begin{align}
\mhY(u):= & \sup_{t\in \RRR^k} <t,u>- \po
\LL_Y(t)-1\pf\label{mhY},\;u\in \RRR^k,\\
\mhYY(x):= & \sup_{t'\in \RRR} t'x- \po \LL_{\mid
Y\mid}(t)-1\pf\label{mhYY},\;x\in \RRR. \end{align} Recall that
the constant $c>0$ appears in assumption (HV) and that
$\Gam_{\mhY}(1/c(z))$ has been defined by (\ref{GamZc}) and
(\ref{mhY}). Our result can be stated as follows.
\begin{theo}\label{T1}
Under assumptions (HV), (HL1) and (HL2), we have almost
surely
\begin{align}\nono
f(z)^{-1}\mathbf{\Delta}_{n,\cc}(z,h_n,\cdot)\;\leadsto\;\Gam_{\mhY}(1/cf(z)).
\end{align}
\end{theo}
A consequence of this result is the following unconsistency
result, for which no proof has been yet provided to the best of
our knowledge: let $K$ be real function on $\RRR^d$ with
\textit{bounded variation and compact support}. The
Nadaraya-Watson regression estimator of $r(z):=\EEE\po Y\mid
Z=z\pf$ is defined as :
$$r_n(z):= \sliin  \frac{K\poo\frac{Z_i-z}{h_n^{1/d}}\pff
}{\sli_{j=1}^n K\poo\frac{Z_j-z}{h_n^{1/d}}\pff}Y_j,\;z\in
\RRR^d.$$ Theorem \ref{T1} entails that, under (HV), (HL1) and
(HL2), the pointwise strong consistency of $r_n$ does
not hold.
\begin{coro}\label{C1}
Under (HV), (HL1) and (HL2), and assuming that $Y\not
\equiv 0$, we have almost surely
$$ \lsn \mid r_n(z)-r(z)\mid_k>0.$$
\end{coro}
\textbf{Proof:} We may assume without loss of generality that $K$
vanishes outside $\Idd$. Consider the random vectors $\wtY_i:=\po
Y_i,1\pf,$ taking values in $\RRR^{k+1}$. Some straightforward
computations show that $\wtY,Z$ satisfy $(HL1)$ and $(HL2)$, and
that, writing $m:=\nabla \LL_{\wtY}(0)$, we have \beq
\mh_{\wtY}(m)=0.\label{ramp}\eeq Moreover, assuming without loss
of generality that $Y$ has a second moment matrix which is
strictly positive, we have $\nabla^2 \LL_{\wtY}>0$ (strictly
positive matrix) on $\RRR^{k+1}$, which ensures that
$$t_0\mapsto \nabla \LL_{\wtY}(t)\Mid_{t=t_0}$$ is a $C^1$
diffeomorphism from $\RRR^{k+1}$ to an open set $O\ni m$. And hence admits an inverse that we write $\nabla \LL_{\wtY}^{-1}$. We
deduce that
$$\mh_{\wtY}(x)=<x,\nabla \LL_{\wtY}^{-1}(x)>-\poo \LL_{\wtY}\po
\nabla\LL_{\wtY}^{-1}(x)\pf-1\pff$$ is continuous in $x$, which implies,
that, for $\e>0$ small enough we have \beq \mathop{\sup_{x\in
\RRR^k}}_{\mmi x-m\mmi_k<\e/f(z)}\;\mid
\mh_{\wtY}(x)\mid_1<\frac{1}{cf(z)}.\label{hape}\eeq
 For $g\in B_{k+1}
(\Idd)$, we shall write $g=(g_k,g_{k+1})$, where $g_{k+1}$ denotes
the last coordinate of $g$ and $g_k\in B_k(\Idd)$ is equal to $g$
without its last coordinate. We shall also write, for a Borel set
$A$ and for $\ell=(\ell_1,\ldots,\ell_k)\in B_k(\Idd)$
\beq\ell(A):=\poo\ili_{\RRR^d}1_{A}d\ell_1,\ldots,\ili_{\RRR^d}1_{A}d\ell_k\pff.\label{integ}\eeq
which is well defined as soon as either $1_A$ or each $g_i$ has bounded variations on
$\RRR^d$. Consider the following mappings:
$$\begin{array}{rcl}
\Psi: \;\poo B_{k+1}(\Idd),\norm_k\pff& \mapsto & \RRR^{k+1}\\
        g&\rar &\pooo\ili_{s\in \Idd} g_k([s,1])dK(s),\ili_{s\in\Idd}g_{k+1}([s,1])
        dK(s)\pff,\\ \\
\Psi': \;\poo B_{k+1}(\Idd),\norm_k\pff & \mapsto & \RRR^{k+1}\\
        g&\rar& \pooo\ili_{s\in\Idd} \ili_{[s,1]}\dot g_k(u) d\lab(u) dK(s),\ili_{s\in \Idd}\ili_{[s,1]}\dot
        g_{k+1}(u)d\lab(u)
        dK(s)\pff, \\ \\
T:\; \RRR^k\times (\RRR-\{0\})&\mapsto &\RRR^k\\
(x_1,\ldots,x_k,x_{k+1})&\rar & \frac{1}{x_{k+1}}\po
x_1,\ldots,x_k).
\end{array}
$$
Also consider
$$\wt{\Gam}_{\mh_{\wt{Y}}}(1/cf(z)):=\aoo g \in
L^2(\Idd)^{k+1},\;\ili_{\Idd}\mh_{\wtY}(g)\le 1/c\aff.$$ We
obviously have $\Psi\poo
f(z)\Gam_{\mh_{\wtY}}(1/cf(z))\pff=\Psi'\poo
f(z)\wt{\Gam}_{\mh_{\wt{Y}}}(1/cf(z))\pff.$ As $K$ has bounded
variations, $\Psi$ is continuous and so is $T\circ \Psi$. Applying
Theorem \ref{T1}, we then deduce that, almost surely
$$r_n(z)=T\circ \Psi \poo
\mathbf{\Delta_{n,\mathfrak{c}}}\pff\;\leadsto \;T\circ\Psi\poo
f(z)\Gam_{\mh_{\wtY}}(1/cf(z))\pff\supset T\circ\Psi'\poo
f(z)\wt{\Gam}_{\mh_{\wt{Y}}}(1/cf(z))\cap
 B_{k+1}(\Idd)\pff.$$ It hence remains to show that
$T\circ\Psi'\poo f(z)\wt{\Gam}_{\mh_{\wt{Y}}}(1/cf(z))\cap
B_{k+1}(\Idd)\pff$ has non empty interior, which shall obviously
imply that, almost surely, $r_n(z)\not \rar r(z)$ as $\nif$. It is
well known that, as $\Psi'$ is continuous, surjective and linear
from the Banach space $\po B_{k+1}(\Idd),\norm_k\pf$ to
$\RRR^{k+1}$, $\Psi'(O)$ is open for every open set $O$. Hence, it
is sufficient to show that
$f(z)\wt{\Gam}_{\mh_{\wt{Y}}}(1/cf(z))\cap
 B_{k+1}(\Idd)$ has nonempty interior in
$\po B_{k+1}(\Idd),\norm_k\pf$. Consider $\e>0$ that appears in
(\ref{hape}). Writing $g_m:\equiv m\in B_{k+1}\po \Idd\pf$, we
have, by (\ref{hape})
$$\mmi f(z)g-f(z)g_m\mmi_k<\e\;\Rightarrow\;\ili_{\Idd}\mh_{\wtY}(g)\le
1/c(z),$$ which concludes the proof. $\Box$  \vk The remainder of
our paper is organised as follows. In \S \ref{Deux}, we introduce
an almost sure approximation of
$\mathbf{\Delta}_{n,\cc}(z,h_n,\cdot)$ by a sum of compound
Poisson processes. This approximation is largely inspired by a
lemma of Deheuvels and Mason \cite{DeheuvelsM6}. We then focus on
these "poissonised" processes and provide some exponential
inequalities on their modulus of continuity. In \S \ref{Trois}, we
establish a Large Deviation Principle (LDP). Then \S \ref{Quatre}
and \S \ref{Cinq} are devoted to proving points $(\ref{i})$ and $(\ref{ii})$
of Theorem \ref{T1} respectively.
\section{A Poisson approximation}\label{Deux}
Recall that $z\in \RRR^d$ is fixed once for all in our problem.
For ease of notation we write \beq \Dnc(z,h,s):= \frac{1}{nh
f(z)}\sliin 1_{[0,s)}\poo \frac{Z_i-z}{h^{1/d}}\pff Y_i,\;\;s\in
\Idd,\;h>0.\label{Dnc}\eeq Throughout this article, we shall refer
to a generic stochastic process $U$, usually called \emph{compound
Poisson process}. It is defined as follows: consider an infinite
i.i.d array $\po \mathfrak{Y}_{ij},\mathfrak{Z}_{ij}\pf_{i\geq
1,\;j\geq 1}$ having the same law as $(Y_1,Z_1)$, as well as a
Poisson random variable with expectation equal to 1 fulfilling
$\eta\;\indep\;\po \mathfrak{Y}_{ij},\mathfrak{Z}_{ij}\pf_{i\geq
1,\;j\geq 1}$ (here $\indep$ denotes stochastic independence). Now
define \beq U(s):=\sli_{j=1}^{\eta}1_{[0,s)}\po
\mathfrak{Z}_{ij}-z\pf\mathfrak{Y}_{ij}.\label{U}\eeq Note that
the law of $U$ is entirely determined by the following property:
\begin{align}
\nono &\text{For each }p\geq 1\text{ and for each partition
}A_1,\ldots,A_p\text{ of }\Idd \text{ we have }:\\
&\EEE\pooo \exp \poo
\sli_{j=1}^p<t_j,U(A_j)>\pff\pfff=\exp\poo\sli_{j=1}^p \PPP(Z-z\in
A_j)\po\LL_{Y\mid A_j}(t_j)-1\pf\pff \label{partoche},\;\po
t_1,\ldots,t_p\pf\in {\po\RRR^k\pf}^p
\end{align}
where $\LL_{Y\mid A_j}(t):=\EEE\poo \exp\po <t,Y>\Mid Z-z\in
A_j\pff$, $j=1,\ldots,p$, $t\in \RRR^k$. Recall the expression
$U(A)$ is understood according to (\ref{integ}). The following
proposition enables to switch the study of the almost behaviour of
the sequence $\po \Dnc(z,h_n,\cdot)\pf\suite$ to that of a
sequence with the following generic term
\begin{align}
\DPnc(h,s):=& \frac{1}{nhf(z)}\sliin U_i(h^{1/d}s),\;s\in
\Idd,\;n\geq 1\label{DPnc},
\end{align}
where the $U_i$ are suitably built independent copies of $U$.
 This result is in the spirit of
Deheuvels and Mason (see \cite{DeheuvelsM6}, Lemma 2.1, or
\cite{DeheuvelsM3}, Proposition 2.1).
\begin{propo}\label{poipoi}
On a probability space rich enough $\wap$ we can construct an
i.i.d. sequence of processes $(U_i)_{i\geq 1}$ having the same law
as $U$ and an sequence $(Y_{i1},Z_{i1})_{i\geq 1}$ having the same
law as $(Y_i,Z_i)_{i\geq 1}$ such that, considering the
$\Dnc(z,h_n,\cdot)$ as built with the sequence
$(Y_{i1},Z_{i1})_{i\geq 1}$ we have almost surely \beq \lsn
\;(nh_n)\; \mmi \Dnc(z,h_n,\cdot)-\DPnc\po h_n,\cdot\pf
\mmi_k<\infty,\label{approximation}\eeq with $\DPnc(\cdot,\cdot)$
defined in (\ref{DPnc}).
\end{propo}
\textbf{Proof}: Denote by $U$ a process having the same law as in
(\ref{U}). Set $\VV_i:=z+h_i^{1/d}\Idd$, $p_i=\PPP\po Z_1\in
\VV_i\pf$, and let $\po Y^{(1)}_{ij}, Z^{(1)}_{ij}\pf_{i\geq
1,\;j\geq 1}$, $\po Y^{(2)}_{ij}, Z^{(2)}_{ij}\pf_{i\geq 1,\;j\geq
1}$, $(\mathbf{b}_i)_{i\geq 1}$, $\po U_i^*(\cdot)\pf_{i\geq 1}$
and $\po v_i\pf_{i\geq 1}$ be families of random elements such
that
\begin{tabbing}
 $ \hskip 5pt $ \= $ (a) $ $ \hskip 5pt $ \= $ \PPP\poo \po Y^{(1)}_{ij},
Z^{(1)}_{ij}\pf \in  B\pff=\PPP\poo (Y_1,Z_1)\in B\Mid Z_1\in \VV_i\pff,\;B\;\text{Borel set},\;i,j\geq 1$. \\
 \nonumber
\> $ (b) $ \> $ \PPP\poo \po Y^{(2)}_{ij},
Z^{(2)}_{ij}\pf \in  B\pff=\PPP\poo (Y_1,Z_1)\in B\Mid Z_1\notin \VV_i\pff,\;B\;\text{Borel set},\;i,j\geq 1$ .\\
\nonumber \> $ (c) $ \> For each $i\geq 1$ we have $\PPP\po v_i=0\pf =1-p_i^{-1}(1-e^{-p_i})$ \\
\>\>and $\PPP\po v_i=k\pf= (k!)^{-1}p_i^{k-1}e^{-p_i},\;k=1,2,\ldots$. \\
\> $ (d) $ \> $ \PPP(\mathbf{b}_i=1)=1-\PPP(\mathbf{b}_i=0)= p_i,\;i\geq 1  $.\\
\nonumber \> $ (e) $ \> The $ U_i^* $ are independent copies of $U$ defined in (\ref{U}) \\
\nonumber\> $ (f) $ \> The union of these five families of random
elements is a stochastically independent family.
\end{tabbing}
In $(e)$, equality in law is understood as an equality with
respect to the $\sig$-algebra $\TT_0$ of $\po
B_k(\Idd),\norm_k\pf$ spawned by the open balls. In $(f)$,
stochastic independence is understood with respect to a suitably
chosen product $\sig$-algebra where each factor is either $\TT_0$,
the Borel $\sig$-algebra of $\RRR^k\times \RRR^d$, or the subsets
of $\{0,1,2,\ldots\}$. First, notice that
$\eta_i^*:=v_i\mathbf{b}_i$ is a Poisson random variable with
expectation $p_i$ for each $i\geq 1$, and that\beq \forall i\geq
1,\;\PPP\po \eta_i^*=\mathbf{b}_i\pf \geq 1-p_i^2
.\label{couplage}\eeq In fact, $\eta_i^*$ and $\mathbf{b}_i$ are a
coupling of a Poisson and Bernouilli random variables
$(\eta,\mathbf{b})$ with expectation $p_i$ such that the
probability $\PPP(\eta=\mathbf{b})$ is maximal. Second, notice
that the following random vectors \beq \po
Y_{ij},Z_{ij}\pf:=1_{\mathbf{b}_i=1}\po Y^{(1)}_{ij},
Z^{(1)}_{ij}\pf+1_{\mathbf{b}_i=0}\po Y^{(2)}_{ij},
Z^{(2)}_{ij}\pf,\;i\geq 1,\;j\geq 1,\label{Yij}\eeq are i.i.d.
with common law equal to $(Y_1,Z_1)$. Moreover, the following
assertions are true with probability one, for each $i\geq 1$: \beq
\forall s\in \Idd,\; 1_{[0,s)}\poo \frac{Z_{ij}-z}{h_i^{1/d}}\pff
Y_{ij}= \sli_{j=1}^{\mathbf{b}_i} 1_{[0,s)}\poo
\frac{Z_{ij}^{(1)}-z}{h_i^{1/d}}\pff Y_{ij}^{(1)}.\label{rd1}\eeq
We now define, for each $i\geq 1$, \beq U_i(s):=U_i^*\poo
[0,s)\cap\{\VV_i-z\}^C\pff+ \sli_{j=1}^{\eta_i^*} 1_{[0,s)}\po
Z_{ij}^{(1)}-z\pf Y_{ij}^{(1)}.\eeq Here, $\VV^C$ denotes the
complement of a given set $\VV\subset \RRR^d$. Some usual
computations on characteristic functions show that the processes
$U_i(\cdot)$ fulfill (\ref{partoche}), and hence are distributed
like $U$. Moreover since $h_{i+p}^{1/d}\Idd \subset h_i^{1/d}\Idd$
for $i\geq 1,\; q\geq 0$, we have almost surely \beq U_i\po
h_{i+q}^{1/d}s\pf =\sli_{j=1}^{\eta_i^*} 1_{[0,s)}\poo
\frac{Z_{ij}^{(1)}-z}{{h_{i+q}}^{1/d}}\pff Y_{ij}^{(1)},\;s\in
\Idd.\label{rd2}\eeq It follows from (\ref{couplage}), (\ref{rd1})
and (\ref{rd2}) that, for each $i\geq 1$,
\begin{align}
 \PPP\pooo U_i\po h_{i+q}^{1/d} \cdot\pf \equiv
1_{[0,\cdot)}\poo\frac{Z_i-z}{h_{i+q}^{1/d}}\pff Y_{i1}\text{ for
each }q\geq 0\pfff \geq \PPP\po \eta_i^*=\mathbf{b}_i\pf
 \geq 1-p_i^2.\label{atob}
\end{align}
Since $p_n =f(z)h_n(1+o(1))$ as $\nif$, and by assumption (HV), we
have $\sli p_i^2<\infty$, which entails, by making use of the
Borel-Cantelli lemma, that (\ref{approximation}) is true with
respect to our construction. $\Box$\vk By Proposition
\ref{poipoi}, proving Theorem \ref{T1} is equivalent to proving a
version of Theorem \ref{T1} with the process $\Dnc(z,h_n,\cdot)$
replaced by their Poisson approxiations $\DPnc(h,\cdot)$. This
will be the aim of \S \ref{Trois}, \S \ref{Quatre} and \S
\ref{Cinq}. In each of these three sections, we shall require the
following exponential inequality for the \textit{absolute
oscillations} of $\DPnc$, which are defined as the oscillations of
the following process : \beq\overline{\DPnc}(h,s):= \frac{1}{nh
f(z)}\sliin \sli_{j=1}^{\eta_i}1_{[0,s)}\poo
\frac{Z_{ij}-z}{h^{1/d}}\pff\mid Y_{ij} \mid_k,\;s\in \Idd,\;n\geq
1\label{VarDPnc}.\eeq Recall that $\mhYY$ has been defined by
(\ref{mhYY}).
\begin{lem}\label{accPoisson}
Given $\dd\in (0,\sqrt{2}-1]$ and $x\geq 0$, there exists
$h_{x}>0$ such that, for each $0<h<h_{x}$ and for each $n\geq 1$,
we have
\begin{align} & \PPP\poooo \mathop{\sup_{s,s'\in\Idd}}_{\mid
s'-s\mid_d\le \dd} \Mid \overline{\DPnc}\po
h,s\pf-\overline{\DPnc}\po h,s'\pf \Mid_1\geq 2d\dd x\pffff \le
{\poo\frac{10}{\dd}\pff}^d\exp\poo -d \dd nh f(z)
\mhYY(x)\pff\label{acc},\\
&\PPP\poo  \overline{\DPnc}(h,1)\geq x\pff\le
\exp\po-nhf(z)\mhYY(x)\pf.\label{total}
\end{align}
\end{lem}
\textbf{Proof:}\lb Given $s$ and $s'\in \Idd$, we write $s\prec
s'$ whenever each coordinate of $s$ is lesser than the
corresponding coordinate of $s'$. Obviously, the
$\overline{\DPnc}\po h,s\pf$ are almost surely increasing in each
coordinate of $s$. First fix $\dd>0$ and set \beq M:=1+\cooo
\frac{3}{(\sqrt{2}-1)\dd}\cfff\label{M}.\eeq We then discretise
$\Idd$ into the following finite grid: \beq
s_{\ii}:=\frac{1}{M}\ii,\;\ii \in
\{0,1,\ldots,M-1\}^d.\label{sii}\eeq By construction, for each $s$
and $s'$ with $\mid s'-s\mid_d\le \dd$, there exists $\ii_s\in
\{0,1,\ldots,M-1\}^d$ fulfilling $s_{\ii_s}\prec s$ and $\mid
s-s_{\ii_{s}}\mid_d\le 1/M$, which entails $\mid s'-s_{\ii_s}\mid
\le 1/M+\dd.$ Hence we can write
\begin{align}
\nono &\PPP\poooo \mathop{\sup_{s,s'\in\Idd}}_{\mid s'-s\mid_d\le
\dd} \Mid
\overline{\DPnc}\po h,s\pf-\overline{\DPnc}\po h,s'\pf \Mid_1\geq 2d\dd x\pffff\\
\nono \le & \PPP\poooo \bculi_{\ii \in \{0,1,\ldots,M-1\}^d} \aoo
\mathop{\sup_{s_{\ii}\prec s',}}_{\mid s'-s_\ii\mid_d\le
\dd+1/M}\Mid \overline{\DPnc}\po h,s'\pf-\overline{\DPnc}\po
h,s_{\ii}\pf
\Mid_1\geq 2d\dd x\aff\pffff\\
\le \nono& M^d \max_{\ii \in \{0,1,\ldots,M-1\}^d}\PPP
\poooo\mathop{\sup_{s_{\ii}\prec s',}}_{\mid s'-s_\ii\mid_d\le
\dd+1/M}\Mid \overline{\DPnc}\po h,s'\pf-\overline{\DPnc}\po
h,s_{\ii}\pf \Mid_k\geq 2d\dd x\pffff.\label{fin}
\end{align}
Now notice that, for each $n\geq 1$, we have\beq
\po\overline{\DPnc}(h,s)\pf_{s\in\Idd}=_{\LL}
\pooo\frac{1}{nhf(z)}\sli_{i=1}^{\eta_n}  1_{[0,s)}\poo
\frac{Z_i-z}{h^{1/d}}\pff \mid Y_i\mid_k\pfff_{s\in
\Idd},\label{rapapou}\eeq where $\eta_n$ is a Poisson random
variable with expectation $n$ and independent of $(Y_i,Z_i)_{i\geq
1}$ (here $=_{\LL}$ stands for the equality in law for processes).
For a Borel set $B\subset \Idd$, write
\begin{align} \nono\overline{\DPnc}(h,B):=&\ili_{\Idd} 1_B(s)d \overline{\DPnc}(h,s)\\
=&\frac{1}{nhf(z)}\sli_{i=1}^{\eta_n}1_B\poo
\frac{Z_i-z}{h^{1/d}}\pff \mid Y_i\mid_k .\label{DPncchA}
\end{align} By the triangle inequality
we have almost surely
\begin{align}
\nono &\mathop{\sup_{s_{\ii}\prec s',}}_{\mid s'-s_\ii\mid_d\le
\dd+1/M}\Mid \ov{\DPnc}\po h,s'\pf-\ov{\DPnc}\po h,s_{\ii}\pf \Mid_1\\
\nono\le &  \mathop{\sup_{s_{\ii}\prec s',}}_{\mid
s'-s_\ii\mid_d\le \dd+1/M} \frac{1}{nhf(z)}\sli_{i=1}^{\eta_n}
 \po 1_{[0,s')}-1_{[0,s_{\ii})}\pf\poo \frac{Z_i-z}{h^{1/d}} \pff \mid Y_i \mid_k \\
 \le & \frac{1}{nhf(z)}\sli_{i=1}^{\eta_n} \po
1_{[0,s_{\ii}^+)}-1_{[0,s_{\ii})}\pf\poo \frac{Z_i-z}{h^{1/d}}
\pff \mid Y_i \mid_k ,\label{rcd}
\end{align}
where $s_{\ii}^+$ is defined by adding $M^{-1}(\co M\dd\cf+2)$ to
each coordinate of $s_i$. Line (\ref{rcd}) is a consequence of the
fact that, if $s_{\ii}\prec s'$ and $\mid s'-s\mid_d\le \dd+1/M$,
then $s_{\ii}\prec s'\prec s_{\ii^+}$. We shall now write
$B_{\ii}:= [0,s_{\ii}^+)-[0,s_{\ii})$. Now choose $t=t(x)$
fulfilling \beq tx-\po \LL_{\mid Y \mid}(t)-1\pf\geq
\frac{1}{2}\mhYY(x).\label{atteint}\eeq By Markov's inequality we
have
\begin{align}
\nono &\PPP \pooo \frac{1}{nhf(z)}\sli_{i=1}^{\eta_n} \po
1_{[0,s_{\ii}^+)}-1_{[0,s_{\ii})}\pf\poo \frac{Z_i-z}{h^{1/d}}
\pff  \mid Y_i
\mid_k\geq 2d \dd x\pfff\\
\nono\le & \exp\po -2d \dd nhf(z)tx \pf \EEE\cooo \exp\pooo
t\sli_{i=1}^{\eta_n} 1_{B_{\ii}}\poo
\frac{Z_i-z}{h^{1/d}}\pff \mid Y_i \mid_k\pfff\cfff\\
\le & \exp\po -2d \dd nhf(z)tx\pf \exp \poo n\po
\LL_{h,\ii}(t)-1\pf\pff,\label{zaza}
\end{align}
with
$$\LL_{h,\ii}(t):=\EEE \cooo \exp \poo t1_{B_{\ii}}\poo \frac{Z_1-z}{h^{1/d}}\pff\mid Y_1
\mid_k\pff\cfff.$$ Note that $(\ref{zaza})$ has been obtained by
conditioning with respect to $\eta_n$. Now, by conditioning with
respect to $E_{\ii,h}:=\{h^{-1/d}(Z_1-z)\in B_{\ii}\}$, and
writing \beq \LL'_{h,\ii}(t):=\EEE\coo \exp \po t\mid Y_1 \mid_k
\pf \Mid E_{\ii,h}\cff,\eeq we obtain
\begin{align}\nono \LL_{h,\ii}(t)-1=&\PPP\po E_{\ii,h}\pf\LL'_{h,\ii}(t)+\poo1-\PPP\po
E_{\ii,h}\pf\pff-1\\
=&\PPP(E_{\ii,h})\poo \LL'_{h,\ii}(t)-1\pff.\label{nunu}
\end{align}
Note that assumptions (HL1) and (HL2) readily entail\beq
\lim_{h\rar 0}\;\max_{\ii \in
\{0,\ldots,M-1\}^d}\Mid\frac{\PPP(E_{\ii,h})\po
\LL_{2,h,\ii}(t)-1\pf}{f(z)\lab(B_{\ii})h\po\LL_{\mid Y
\mid_k}(t)-1\pf}-1 \Mid=0.\label{fundam}\eeq Choose $h_{x}>0$
small enough so that the quantity involved in (\ref{fundam}) is
lesser that $\sqrt{2}-1$ and notice that for each $\ii$ we have
$\lab(B_\ii)\le d(\dd+1/M)\le \sqrt{2}d\dd$ by (\ref{M}). By
combining (\ref{rcd}), (\ref{zaza}) and (\ref{nunu}), we conclude
that, for all $0<h<h_{x}$,
\begin{align}
\nono &\max_{\ii \in \{0,1,\ldots,M-1\}^d}\PPP
\poooo\mathop{\sup_{s_{\ii}\prec s',}}_{\mid s'-s_\ii\mid_d\le
\dd+1/M}\Mid \ov{\DPnc}\po h,s'\pf-\ov{\DPnc}\po h,s_{\ii}\pf
\Mid_k\geq
2d\dd x\pffff\\
\le & \exp\poo -2d \dd nhf(z) \po tx-\LL_{\mid
Y\mid_k}(t)+1\pf\pff,
\end{align}
whence, by (\ref{atteint}) we get
\begin{align}
\nono &\PPP\poo \mathop{\sup_{s,s'\in\Idd}}_{\mid s'-s\mid_d\le
\dd} \Mid
\ov{\DPnc}\po h,s\pf-\ov{\DPnc}\po h,s'\pf \Mid_k\geq 2d\dd x\pff\\
\nono \le & M^d\exp\poo -d \dd nhf(z)  \mhYY(x)\pff\\
\nono\le& \poo1+\frac{3}{(\sqrt{2}-1)\dd}\pff^d\exp\poo -d \dd nh
f(z)
\mhYY(x)\pff\\
\nono \le &{\poo\frac{10}{\dd}\pff}^d\exp\poo -d  \dd nh f(z)
\mhYY(x)\pff.
\end{align}
This concludes the proof of Lemma \ref{accPoisson}. $\Box$
\section{Large deviations for $\DPnc(h_n,\cdot)$}\label{Trois}
In this section, we establish a Large Deviation Principle (LDP)
for the sequence of processes $\DPnc(h_n,\cdot)$. For the
definition of large deviations for sequences for bounded
stochastic processes and of a (good) rate function, we refer to
Arcones \cite{Arcones1}.
\subsection{Some tools in large deviation theory}
We begin this subsection with some well known properties (see,
e.g., \cite{Deheuvels4}, Lemma 2.1, or Borovkov \cite{Borovkov1}
just above the main Theorem) of $\mhY$ and $\mhYY$ given in
(\ref{mhY}) and (\ref{mhYY}) respectively.
\begin{fact}\label{propmh}
The functions $\mhY$ and $\mhYY$ are positive convex. Moreover,
since $\LL_Y$ is finite on $\RRR^k$, we have
\begin{align}
\nono \lim_{\mid u \mid_k\rar \infty} \frac{\mhY(u)}{\mid u\mid
_k}=&\infty,\\
\nono \lim_{\mid x \mid\rar \infty} \frac{\mh_{\mid Y_k
\mid}(x)}{\mid x\mid }=&\infty.
\end{align}
\end{fact}
Arcones (see \cite{Arcones1}, Theorem 3.1) has established a very
useful criterion to establish a LDP for processes in $B_k(\Idd)$
(actually only with $k=1$ but the extension of his results to
$k>1$ is straightforward). We cannot make a direct use of his
Theorem 3.1 and shall make use of a slight modification of it. To
state this modification, we shall introduce some more notations.
For each integer $p\geq 1$, consider a finite grid
\begin{align}
\nono S_p=&\aoo s_{\jj,p},\jj\in \ppd\aff\\
:=&\aoo 2^{-p}\jj,\;\jj\in\ppd\aff. \end{align} and consider its
associated partition of $\Idd$ into hypercubes, namely
\begin{align}
C_{\jj,p}:= \coo 2^{-p}(\jj-1),2^{-p}\jj\pff,\;\jj\in
\ppd\label{Cjjp}.
\end{align}
Here we have written $\jj-1=(j_1-1,\ldots,j_d-1)$.
 Now for each integer $p\geq 1$ and for each $g\in B_k(\Idd)$
 write
 \beq g^{(p)}(s):= g(s_{\jj,p}),\;s\in C_{\jj,p},\;\jj\in\ppd.\label{gp}\eeq
 The following proposition is a straightforward variation of
 Theorem 1 of Arcones \cite{Arcones1}, and is written according to the
 notation of that theorem (in particular, we refer to
 \cite{Arcones1} for a definition of the outer probability
 $\PPP^*$).
\begin{propo}\label{propo: ecrit}
Let $(X_n)\suite$ be a sequence of stochastic processes and let
$(\e_n)\suite$ be a sequence of constants fulfilling
$\e_n>0,\;n\geq 1$ and $\e_n\rar 0$ as $\nif$. Assume that the
following conditions are satisfied.
\begin{enumerate}
\item  The sequence of stochastic processes $(X_n^{(p)})\suite$
satisfies the LDP for $(\e_n^{-1})\suite$ and for a rate function
$J_p$ on $\po B_k(\Idd),\norm_k\pf$. \item For each $\tau>0$ and
$M>0$ there exists an integer $p\geq 1$ satisfying
\[\lsn \e_n\log\poo\PPP^*\poo\max_{\jj\in\ppd}\sup_{s\in C_{\jj,p}}\mid X_n(s)-X_n(s_{\jj,p})\mid_k \geq
\tau\pff\pff\le -M.\]
\end{enumerate}
Then $(X_n)\suite$ satisfies the LDP for $(\e_n^{-1})\suite$ and
for the following rate function.
$$J(g):=\sup_{p\geq
1}J_p\poo g^{(p)}\pff,\;g\in B_k(\Idd).$$
\end{propo}
\textbf{Proof:} The proof follows the same lines as in the proof
of Theorem 3.1 of Arcones \cite{Arcones1}. We omit details for
sake of briefness. $\Box$\vk For $g=(g_1,\ldots,g_k)\in B_k(\Idd)$
and $A$ Borel set, we shall write \beq g(A):=\pooo\ili_{\Idd}
1_A(s)dg_1(s),\ldots,\ili_{\Idd} 1_A(s)dg_k(s)\pfff,\label{gA}\eeq
which is valid as long as $1_A$ or each $g_l$ have bounded
variations. We shall now consider the following (rate) functions
on $\po B_k(\Idd),\norm_k\pf$ that will play the role of
successive approximations of $\JhY$: given $p\geq 1$ and $g\in
B_k(\Idd)$ we set \beq \JhY^{(p)}(g):=\sli_{\jj\in\ppd}
\lab(C_{\jj,p})\mhY\poo
\lab(C_{\jj,p})^{-1}g(C_{\jj,p})\pff.\label{JhYp}\eeq The
following fact is a straightforward extension to the multivariate
case of Proposition 2.1 in \cite{Varron3}. Recall that $\JhY$ has
been defined through (\ref{JZ}) and (\ref{mhY}).
\begin{fact}\label{lscBVk}
For any $g\in B_k(\Idd)$ we have \beq \JhY(g)=\lim_{p\rar\infty}
\JhY^{(p)}(g).\eeq As a consequence, $\JhY$ is lower
semicontinuous on $B_k(\Idd)$.
\end{fact}
Our next lemma states that the function $\JhY$ (recall (\ref{JZ}))
is a "rate" function.
\begin{lem}\label{rate}The sets $\Gam_{\JhY}(a),\;a\geq 0$ are compact subsets of $\po B_k(\Idd),\norm_k\pf$. In other words,
$J_{\mhY}$ is a rate function in $\po B_k(\Idd),\norm_k\pf$.
\end{lem}
\textbf{Proof :} By Fact \ref{propmh} we have $\mid x \mid_k\le
\mid x \mid_k1_{\mid x\mid_k\le M}\wedge \JhY$ for some $M>0$ and
for each $x$. Hence, for any $g\in \Gam_{\JhY}(a)$ we have (recall
that $\lab$ stands for the Lebesgue measure)
\begin{align}
\ili_{\Idd} \mid  g'\mid_kd\lab=&\ili_{\mid g'\mid_k \le M}\mid
g'\mid_kd\lab +\ili_{\mid g'\mid_k > M}\JhY\po g'\pf d\lab\\
\le& M+a,\label{variation}
\end{align}
from where we conclude that $\Gam_{\JhY}$ is relatively compact in
$B_k(\Idd)$. It is also closed in $B_k(\Idd)$ by  a combination of
Fact \ref{lscBVk} and (\ref{variation}), which proves Lemma
\ref{rate}.$\Box$\lb
\subsection{A large deviation principle }
In this subsection, we state and prove a large deviation principle
that will play a crucial role in the sequel of our proof of
Theorem \ref{T1}. This LDP is stated as follows:
\begin{propo}\label{ldp}
Under assumptions $(HV),\;(HL1)-(HL2)$, the sequence $\poo
\DPnc(h_n,\cdot)\pff\suite$ satisfies the LDP in $B_k(\Idd)$ for
$(\e_n^{-1})\suite=\po (n h_nf(z))^{-1}\pf\suite$ and for the rate
function $\JhY$.
\end{propo}
\textbf{Proof:} As we shall make use of Proposition \ref{propo:
ecrit}, we have to check conditions 1 and 2 of that proposition,
which will be the aim of the following lemmas. Notice that, almost
surely, we have
\begin{align}
\nono \DPnc(h,C):=&\frac{1}{nhf(z)}\sliin
\sli_{j=1}^{\eta_i}1_C\poo \frac{Z_{i,j}-z}{h^{1/d}}\pff
Y_{i,j}\label{DPnchA}, \;C\text{ Borel },\end{align} with
$\DPnc(h,C)$ defined according to (\ref{integ}). Our proof is
divided in two steps, where we shall respectively verify
conditions 1 and 2 of Proposition \ref{propo: ecrit}.\vk
\textbf{Step 1}: To check condition 2 of Proposition \ref{propo:
ecrit}, we shall make use of Lemma \ref{accPoisson}, which readily
entails, for fixed $p\geq 1$ and $\tau>0$, and for all $n\geq
n(p,\tau)$:
\begin{align}
  \nono &\PPP\poo\max_{\jj\in\ppd}\sup_{s\in C_{\jj,p}}\mid \DPnc(h_n,s)-\DPnc(h_n,s_{\jj,p})\mid_k \geq
\tau\pff\\
\nono \le & 10^d2^{pd}\exp\poo-d2^{-p}nh_nf(z)\mh_{\mid
Y\mid_k}\po d^{-1}2^{p-1}\tau\pf\pff.
\end{align}
Now fix $M>0$ and $\tau>0$. By Fact (\ref{propmh}), we have, for
all large $p$:
$$\frac{\mh_{\mid
Y\mid_k}\po d^{-1}2^{p-1}\tau\pf}{d^{-1}2^{p-1}\tau}>4M\tau,$$
which implies that condition 2 of Proposition \ref{propo: ecrit}
is verified.\lb \textbf{Step 2}: To check condition 1 of
Proposition \ref{propo: ecrit}, we shall require the following
preliminary lemma.
\begin{lem}\label{fidi1}
For any sequence $(h_n)\suite$ fulfilling $h_n\rar 0$ and
$nh_n\rar \infty$, and for any fixed $p\geq 1$, the sequence of
random vectors of ${\RRR^{k}}^{2^p}$ \beq  \poo
\DPnc(h_n,C_{\jj}),\;\jj \in \ppd \pff_{n\geq 1}\label{gro}\eeq
satisfies the LDP for the sequence
$(\e_n^{-1})\suite:=((nh_nf(z))^{-1})\suite$ and the following
rate function $$\begin{array}{rcl}
  J^{(p)}: {(\RRR^k)}^{2^p}& \mapsto & \RRR \\
  x& \mapsto & \sli_{\jj \in \ppd}^p \lab(C_{\jj,p})\mhY\poo
\lab(C_{\jj,p})^{-1}x_\jj\pff. \\
\end{array}$$ Here we write
$x=(x_{\jj},\;\jj\in \ppd)$, with $x_{\jj}\in {\RRR^{k}}$ for each
$\jj \in \ppd$.
\end{lem}
\textbf{Proof:} The proof of Lemma \ref{fidi1} is divided into
three steps. The two first steps deal with a single component of
the random vectors written in (\ref{gro}).\lb \textit{Step 1} :In
our \textit{first step}, we make an additional assumption on
$\LL_Y$, which allows us to make a full use of the Gärtner-Ellis
theorem (see, e.g., \cite{DemboZ}, p. 44). $$ (H_0):\;\;\forall
x\in\RRR^k\text{ fulfilling }\mhY(x)<\infty,\;\exists \eta \in
\RRR^k,\;x=\nabla \LL_Y (\eta).$$
\begin{lem}\label{fidi11}
Assume that $(H_0)$ is true in addition to the assumptions of
Theorem \ref{T1}. Then, for each $p\geq 1$ and $\jj \in \ppd$, the
sequence
$$\poo \DPnc(h_n,C_{\jj,p})\pff\suite $$
satisfies the LDP for the sequence $(nh_n f(z))^{-1}$ and the rate
function $\lab(C_{\jj,p})\mhY\po \lab(C_{\jj,p})^{-1}\cdot\pf.$
\end{lem}
\textbf{Proof of Lemma \ref{fidi11}:}\lb We shall first show that,
for each $t\in \RRR^k$, we have \beq \limn \frac{1}{nh_nf(z)}\log
\pooo \EEE\poo \exp
<t,nh_nf(z)\DPnc(C_{\jj,p},h_n)>\pff\pfff=\lab\po C_{\jj,p}\pf
\poo \LL_Y\poo \lab\po C_{\jj,p}\pf^{-1}
t\pff\pff.\label{moments}\eeq To show this, we start from the
equality (\ref{rapapou}) to obtain by convolution:
\begin{align}
\nono \log \pooo \EEE\poo \exp
<t,nh_nf(z)\DPnc(C_{\jj,p},h_n)>\pff\pfff=n\log \pooo \EEE\poo
\exp <t,nh_nf(z)U\po h_n^{1/d}C_{\jj,p}\pf>\pff\pfff.
\end{align} Recall that $U$ has been defined in (\ref{U}). Next, we use the characterisation (\ref{partoche}), which is
applied to the simple partition $\poo
h_n^{1/d}C_{\jj,p},\;\Idd-h_n^{1/d}C_{\jj,p}\pff$. Using that
relation with $t_1=t$ and $t_2=0$, we obtain
\begin{align}
\nono &\log \pooo \EEE\poo \exp
<t,nh_nf(z)\DPnc(C_{\jj,p},h_n)>\pff\pfff\\
\nono=&n\PPP\poo Z-z\in h_n^{1/d}C_{\jj,p}\pff\EEE\poo \exp\po
<t,Y>\pf\Mid Z-z\in h_n\in z+h_n^{1/d}C_{\jj,p}\pf\pff-1\pff.
\end{align}
Hence (\ref{moments}) follows from assumptions $(HL1)-(HL2)$.\lb B
Lemma 2.3.9 in \cite{DemboZ}, p 46, we know that $(H_0)$ implies
that the set of \textit{exposed points} of $\mhY$ is equal to
$\{x\in \RRR^k,\; \mh(x)<\infty\},$ from where the proof of Lemma
\ref{fidi11} is concluded by an application of the Gärtner-Ellis
theorem (see, e.g., \cite{DemboZ}, p. 44).\nocite{Ellis} $\Box$
\lb \textit{Step 2}: In our \textit{second step}, we shall get rid
of assumption $(H_0)$, which is unfortunately not verified in all
situations (for example, take $k=1$, $Y\equiv 1$, which leads to
$\LL_Y(t)=\exp(t)$, $t\in \RRR$ and $\mhY(0)=1$, but $(H_0)$ is
not satisfied for $x=0$).
\begin{lem}\label{fidi12}
Lemma \ref{fidi11} is true without making assumption $(H_0)$.
\end{lem}
\textbf{Proof of Lemma \ref{fidi12}}: First notice that the
"closed sets" part of the LDP stated in Lemma \ref{fidi11} can be
proved by making use of the Gärtner-Ellis theorem, without making
assumption $(H_0)$. Only the "open sets" part of Lemma
\ref{fidi11} needs assumption $(H_0)$, since it implies that the
set of \textit{exposed points} of $\mhY$ is equal to $\{x\in
\RRR^k,\; \mhY(x)<\infty\}$. We only need to prove that, without
assumption $(H_0)$, for any open set $O\subset \RRR^k$ with
$\mhY(O)<\infty$ (nontrivial case), we have \beq \lin
\frac{1}{nh_nf(z)}\log \pooo \PPP\poo \DPnc\po C_{\jj,p},h_n\pf\in
O\pff\pff\geq -\mhY(O).\label{ouvy}\eeq To achieve this goal, we
shall slightly modify the $Y_{i,j}$ by adding small Gaussian
random vectors. Fix $O\subset \RRR^k$ open, with $\mhY(O)<\infty$,
and $\dd>0$. There exists $x\in O$ and $\dd_1\in (0,\dd)$ such
that $B(x,2\dd_1)\subset O$ and $\mhY(O)\le \mhY(B(x,2\dd_1))\le
\mhY(x)\le \mhY(O)+\dd<\dd_1^{-1}$. Here $B(x,\e)$ denotes the
open ball with centre $x$ and radius $\e$. Now introduce an array
$\po \zeta_{ij}\pf_{i,j\geq 1}$ of $\RRR^k$ valued standard random
vectors, that are independent of the array $\po
Y_{i,j},Z_{i,j}\pf_{i,j\in \NNN}$. Also define
\begin{align}
\nono \DPnc'\po C_{\jj,p},h_n\pf:=&\frac{1}{nh_n f(z)}\sliin
\sli_{j=1}^{\eta_i} 1_{C_{\jj,p}}\poo \frac{
Z_{ij}-z}{h_n^{1/d}}\pff 
\zeta_{ij},\\
\nono \DPnc''\po C_{\jj,p},h_n\pf:=&\frac{1}{nh_n f(z)}\sliin
\sli_{j=1}^{\eta_i} 1_{C_{\jj,p}}\poo
\frac{Z_{ij}-z}{h_n^{1/d}}\pff \po
Y_{ij}+\dd_1^2\zeta_{ij}\pf\\
\nono =& \DPnc\po C_{\jj,p},h_n\pf+\dd_1^2\DPnc'\po
C_{\jj,p},h_n\pf.
\end{align}
We shall first show that the vector $Y+\zeta$ fulfills assumptions
$(H_0)$. To prove this first notice that
$\LL_{Y+\zeta}=\LL_{Y}\LL_{\zeta}$, which holds since $Y$ and
$\zeta$ are independent conditionally to $Z$. Obviously we have,
since $\zeta\indep Z$,
$$\LL_{\zeta}(t)=\exp \poo \frac{1}{2}\mid t\mid_k^2\pff,$$
which shows that $\zeta$ fulfills $(H_0)$. Moreover, by Jensen's
inequality we have
$$\LL_Y(t)\geq \exp\po <m_Y,t>\pf ,\;t\in \RRR^k,$$
where $m_Y=\EEE\po Y\Mid Z=z\pf$, which leads to \beq
\LL_{Y+\dd_1^2\zeta}(t)\geq \exp\poo <m_Y,t>+\frac{\dd_1^4}{2}\mid
t\mid_k^2\pff.\label{masa}\eeq
 Now consider
$x\in\RRR^k$, and define the function $g(t)=<x,t>-\po
\LL_{Y+\dd_1^2\zeta}(t) -1\pf.$ By (\ref{masa}) we have $g(t)\rar
-\infty$ as $\mid t\mid_k\rar \infty$. Hence, the continuous and
differentiable function $g$ admits a maximum at some $\eta \in
\RRR^k$ fulfilling $0=\nabla g(\eta)=y-\nabla \LL_{Y+\zeta}(x)$.
This proves that the vector$\zeta$ fulfills $(H_0)$ and hence, by
Lemma \ref{fidi11} we have:
\begin{align} \lin \frac{1}{nh_n f(z)}\log \pooo \PPP \poo
\DPnc''\po
C_{\jj,p},h_n\pf\in O\pff\pfff\geq &-\mh_{Y+\dd_1^2\zeta}(O),\\
\lsn \frac{1}{nh_n f(z)}\log \pooo \PPP \poo \Mmi\DPnc'\po
C_{\jj,p},h_n\pf\Mmi_k\geq \dd_1^{-1}\pff\pfff\le &-\inf_{\mmi x
\mmi_k \geq \dd_1^{-1}}\mh_{\zeta}(x)\\
\nono\le &-\dd_1^{-1}.
\end{align}
The last inequality holds for $\dd_1>0$ small enough, by Fact
\ref{propmh}, replacing $Y$ by $\zeta$. Hence, by the triangle
inequality, we have for all large $n$:
\begin{align}
 \nono \PPP\poo \DPnc(C_{\jj,p},h_n) \in O\pff\geq& \PPP\poo
 \Mmi\DPnc(C_{\jj,p},h_n)-x\Mmi_k<2\dd_1\pff\\
 \nono \geq &\PPP\poo \poo
 \Mmi\DPnc''(C_{\jj,p},h_n)-x\Mmi_k<\dd_1\pff-\PPP\poo
 \dd_1^2\Mmi\DPnc'(C_{\jj,p},h_n)\Mmi_k>\dd_1\pff\\
 \nono \geq & \exp \pooo -nh_n f(z)\poo \dd+\mh_{Y+\zeta}\po
 B(x,\dd_1)\pf\pff\pfff-\exp\poo -nh_nf(z)\dd_1^{-1}\pff\\
 \geq & \exp \pooo -nh_n f(z)\poo \dd+\mh_Y(x)\pff\pfff-\exp\poo
 -nh_nf(z)\dd_1^{-1}\pff\label{explain1}\\
 \geq & \frac{1}{2}\exp \pooo -nh_n f(z)\poo
 2\dd+\mh_Y(O)\pff\pfff\label{explain2}.
\end{align}
Note that (\ref{explain1}) is a consequence $\mh_{Y+\zeta}\le
\mhY$, which follows directly from $\LL_{\zeta}\geq 1$. Also,
(\ref{explain2}) is a consequence of $\mhY(x)\le \mhY(O)+\dd$
together with $\dd_1^{-1}>\mhY(O)+2\dd$, which is true by the
choice of $\dd_1$. The proof of Lemma \ref{fidi12} is then
concluded since $O$ and $\dd$ are arbitrary.$\Box$\vk
\noindent\textit{Step 3}: The proof of Lemma \ref{fidi1} by a
tensorisation argument brought by Lynch an Sethuraman. Since, for
each $n$, the collection
$$\DPnc(h_n,C_{\jj}),\;\jj \in \ppd$$ is independent, and since
each sequence $\poo \DPnc(h_n,C_{\jj})\pff\suite$ satisfies the
LDP with the rate function $\lab(C_{\jj,p})\mhY\po
\lab(C_{\jj,p})^{-1}\cdot\pf$. Then Lemma \ref{fidi1} is proved by
applying Lemma 2.8 in \cite{Lynch}. $\Box$\vk A direct consequence
of Lemma \ref{fidi1} is that condition 1 of Proposition
\ref{propo: ecrit} is satisfied, as shows our next lemma.
\begin{lem}\label{fidi2}
If $h_n\rar 0$ and $nh_n\rar \infty$, then the sequence of
processes
$$\poo\DPnc^{(p)}(h_n,\cdot)\pff\suite$$
satisfies the LDP for $\e_n:=(nh_nf(z))^{-1}.$ and for the rate
function $\JhY^{(p)}$.
\end{lem}
\textbf{Proof:} The proof is a straightforward application of the
contraction principle (see, e.g., \cite{Arcones1}, Theorem 2.1),
considering, for fixed $p$, the following application, from
${\RRR^k}^{2^{pd}}$ to $\po B_k(\Idd);\norm_k\pf$ (here we write
$x=(x_{\jj},\;\jj \in \ppd)$, with each $x_\jj$ belonging to
$\RRR^k$)
$$\begin{array}{rcl} \mcal{R}_p(x):\Idd&\mapsto& [0,\infty)\\
s&\rar&\sli_{C_{\jj,p}\subset [0,s_1]\times\ldots\times
[0,s_p]}x_{\ii}.\Box
\end{array}$$
We conclude the proof of Proposition \ref{ldp} by combining Step 1
and Step 2 with Proposition \ref{propo: ecrit}.$\Box$
\section{Proof of point $(\ref{i})$ of Theorem \ref{T1}}\label{Quatre}
We shall make use of some usual blocking arguments along the
following subsequence: \beq n_k:= \cooo\exp\poo
k\exp\po-\sqrt{\log k}\pf\pff\cfff,\eeq with associated blocks
$N_k:=\{\nkm+1,\ldots,\nk\}$. Here, $[u]$ denotes the only integer
fulfilling $[u]\le u \le[u]+1.$ We point out two key properties of
$(n_k)_{k\geq 1}$: \beq \limk \frac{\nk}{\nkm}=1,\;\;\limk
\frac{\log_2 n_k}{\log k}=1\label{propnk}.\eeq For any $\e>0$ and
$A\subset B_k(\Idd)$, we shall write: \beq A^\e:=\aoo g\in
B_k(\Idd),\; \inf_{g'\in A}\mmi
g-g'\mmi_k<\e\aff.\label{Aepsilon}\eeq Now, recalling the
definition of $\DPnc$ in (\ref{DPnc}), we define the following
normalised Poisson processes that will play a crucial role in our
blocking arguments. \beq \HH_n(s):=\frac{1}{n_k\hk
f(z)}\sli_{i=1}^n U_i(\hk s),\;k\geq 1,\;n\in N_k,\;s\in
\Idd\label{HHn}.\eeq Fix $\e>0$. We shall proceed in two steps:
first, we will prove that, we have almost surely, ultimately as
$\nif$, \beq \HH_n\in {\GhY(1/cf(z))}^{2\e},\label{l1}\eeq then we
shall show that almost surely: \beq \limk \max_{n\in N_k} \mmi
\HH_n(\cdot)-\DPnkc(\hk,\cdot)\mmi_k\le 3\e \label{l2}.\eeq
\textbf{Step 1}: We first prove (\ref{l1}). In order to make use
of usual blocking arguments along the blocks $N_k$ we shall first
show that \beq \limk \max_{n\in N_k}\PPP\poo \mmi
\HH_n(\cdot)-\DPnkc(\hk,\cdot)\mmi_k>\e\pff=0.\label{block1}\eeq
To prove this, choose $k\geq 1$ and $n\in N_k$ arbitrarily. A
rough upper bound gives (excluding the trivial case where
$n=\nk$).
\begin{align}
\nono \PPP_{n,1}:=&\PPP\poo
\mmi\HH_n(\cdot)-\DPnkc(\hk,\cdot)\mmi_k>\e\pff\\
\nono\le & \PPP\poo \sli_{i=1}^{\nk-n}\sli_{j=1}^{\eta_i}\mid
Y_{i,j}\mid_k>\e\frac{\nk}{\nk-n}(\nk-n)\hk f(z)\pff\\
=&\PPP\poo
\overline{\Delta\Pi}_{\nk-n,\cc}(\hk,\cdot)>\e\frac{\nk}{\nk-n}\pff
\end{align}
Now making use of point (\ref{total}) of Proposition
\ref{accPoisson} with $x:=\e\nk/(\nk-n)$ we get, for all large $k$
and for each $n\in N_k$ with $n\neq \nk$, \beq \PPP_n\le \exp\poo
-\e\nk \hk f(z)
\frac{\nk-n}{\e\nk}\mhYY\po\frac{\e\nk}{\nk-n}\pf\pff.\eeq Now, as
$\nk-n\geq \nk-\nkm$, $\nk/(\nk-\nkm)\rar \infty$ and by Fact
\ref{propmh} we readily infer (\ref{block1}).\lb We are now able
to make use of a well known maximal inequality (see, e.g.,
Deheuvels and Mason \cite{DeheuvelsM2}, Lemma 3.4) to conclude
that, for all large $k$,
\begin{align}
\nono\PPP_{k,2}:=&\PPP\pooo \bculi_{n\in N_k} \HH_n\notin
\GhY^{2\e}\pff\\
\nono\le& 2 \PPP\poo \HH_{\nk}\notin \GhY^{\e}\pff\\
=&2 \PPP\poo \DPnkc(\hk,\cdot)\notin \GhY^{\e}\pff.
\end{align}
Applying proposition \ref{ldp} to the closed set
$F:=B_k(\Idd)-{\po\GhY\pf}^{\e}$, which satisfies $\JhY(F)\geq
(1+3\alp)/cf(z)$ for some $\alp>0$ (by lower semi continuity of
$\JhY$) we get, ultimately as $\kif$,
\begin{align}
\nono \PPP_{k,2}\le& 2\exp\poo - \frac{\nk\hk(1+2\alp)}{cf(z)}\pff\\
\le &\exp\po -(1+\alp)\log\log \nk\pf,\label{sophax}
\end{align}
where (\ref{sophax}) is a consequence of assumption $(HV)$. By
(\ref{propnk}), we conclude that $(\PPP_{k,2})_{k\geq 1}$ is
summable, which proves (\ref{l1}) by making use of the
Borel-Cantelli lemma.\lb  \textbf{Step 2} To prove (\ref{l2}) we
shall make use of the following almost sure equality \beq
\DPnc(h_n,s):=\frac{\nk\hk}{nh_n}\HH_n\po \frac{\hk}{h_n}s\pf.\eeq
By (\ref{propnk}) together with $(HV)$ we straightforwardly infer
that \beq \limk \max_{n\in N_k} \Mid
\frac{\nk\hk}{nh_n}-1\Mid=0,\;\limk \max_{n\in N_k}
\frac{\hk}{h_n}=1.\label{oh1}\eeq Moreover, making use of
(\ref{variation}), we infer that \beq \lim_{T\rar 1,\;\rho\uparrow
1}\sup_{g\in \GhY}\;\mmi
Tg(\rho\cdot)-g(\cdot)\mmi_k=0\label{oh2}.\eeq Hence, (\ref{l2})
follows from a combination of (\ref{oh1}), (\ref{oh2}) and
(\ref{l1}) together with the triangle inequality.\lb The proof of
point $(\ref{i})$ of Theorem \ref{T1} is concluded by combining
(\ref{l1}) and (\ref{l2}) and recalling that $\e>0$ was
arbitrary.$\Box$
\section{Proof of point $(\ref{ii})$ of Theorem \ref{T1}}\label{Cinq}
\newcommand{\hhnk}{h_{\overline{n}_k}}
We introduce the following subsequence $$\ovnk:=k^{2k},\;k \geq
1.$$ Obviously, $\ovnk$ satisfies the following properties: \beq
\log_2 \ovnk=\log k+\log_2 k+\log
2,\;\;\ovnk/\ovnkm=e^{-2}k^{-2}(1+o(1))\label{propovnk}.\eeq we
also shall write $v_k:= \ovnk-\ovnkm$. Now define the sequence
$$\HH_k':=\frac{1}{v_k\hhnk f(z)}\sli_{\ovnkm+1}^{\ovnk}1_{[0,\cdot]}\poo
\frac{Z_i-z}{\hhnk^{1/d}}\pff Y_i.$$ Now choose $\e>0$ and $g\in
\Gam_{\mhY}(1/cf(z))$ arbitrarily. We shall prove that, with
probability one \beq \lsn \Mmi \HH_k'-g\Mmi\le 2\e,
\label{elouan}\eeq which would conclude the proof of point $(\ref{ii})$
of Theorem \ref{T1} by a classical compactness argument. Obviously
$g$ satisfies \beq \lim_{\rho\rar 1}\mmi g(\rho
\cdot)-g(\cdot)\mmi_k=0\label{z1}.\eeq Some routine analysis also
shows that, for some $\alp>0$ we have $J\po
g^\e\pf<(1-2\alp)/cf(z).$ By (\ref{propovnk}) we have
$v_k\hhnk\rar\infty$ as $\kif$. Hence, by Proposition (\ref{ldp}),
which we apply to the open ball $g^{\e}$ we obtain, for all large
$k$ \begin{align} \nono \PPP\poo \HH'_k\in g^{\e}\pff\geq
&\exp\poo -\frac{v_k\hhnk f(z)(1-2\alp)}{c f(z)}\pff\\
\nono \geq & \exp\poo -\log k+\log_2 k+\log 2\pff,
\end{align}
where the last inequality is a consequence of (\ref{propovnk}). As
the $(\HH_k')_{k\geq 1}$ are independent, the Borel-Cantelli lemma
entails, almost surely,\beq \mmi \HH_k'-g\mmi_k\le \e\text{ for
all large }k.\label{asasa}\eeq
To conclude the proof, notice that
\begin{align} \nono\HH_k'=&\frac{v_k}{\ovnk}\HH_k'+\frac{1}{\ovnk \hhnk
f(z)}\sli_{i=1}^{\ovnkm}1_{[0,\cdot]}\poo
\frac{Z_i-z}{\hhnk^{1/d}}\pff Y_i\\
=:&\frac{v_k}{\ovnk}\HH_k'+\zeta_k.
\end{align}
Hence, if we show that $\mmi \zeta_k\mmi_k \rar 0$ almost surely,
then (\ref{elouan}) will follow by noticing that $v_k/\ovnk\rar 1$
and applying both (\ref{z1}) and (\ref{asasa}). Noticing that
$$\mmi \zeta_k\mmi_k\le
\frac{\ovnkm}{\ovnk}\overline{\Delta\Pi}_{\ovnkm,\mathfrak{c}}(\hhnk,1),$$
we readily infer, by (\ref{propovnk}) and point (\ref{total}) of
Lemma \ref{accPoisson}, that
$\PPP\po\mmi\zeta_k\mmi_k>\dd\pf=O(k^{-2})$ for any $\dd>0$. This
concludes the proof of point $(\ref{ii})$ of Theorem \ref{T1}.

\end{document}